\documentclass[preprint]{elsarticle}

\usepackage{amssymb}
\usepackage{amsfonts}
\usepackage{mathrsfs}
\usepackage{algorithm}
\usepackage{algorithmic}
\usepackage{multirow}
\usepackage{xcolor}
\usepackage{txfonts}

\usepackage{graphicx}
\usepackage{subfigure}

\newtheorem{theorem}{Theorem}
\newtheorem{lemma}{Lemma}
\newtheorem{proposition}{Proposition}
\newtheorem{remark}{Remark}

\newtheorem{assumption}{Assumption}

\newcommand{\argmin}{\mathop{{\rm arg}\min}}

\journal{}

\begin{document}

\begin{frontmatter}

\title{Error analysis of  regularized least-square regression with Fredholm kernel}

\author[hubu]{Yanfang Tao}
\ead{tyf3122@163.com}
\author[hzau1]{Peipei Yuan}
\ead{yuanpeipei@webmail.hzau.edu.cn}

\author[hzau2]{Biqin Song\corref{cor}}
\ead{biqin.song@mail.hzau.edu.cn}

\cortext[cor]{Corresponding author.}

\address[hubu]{Faculty of Mathematics and Statistics, Hubei University, Wuhan 430062, China}
\address[hzau1]{College of Engineering,  Huazhong Agricultural University,    Wuhan 430070,
China}
\address[hzau2]{College of Science,  Huazhong Agricultural University,    Wuhan 430070,
China}

\begin{abstract}
Learning with Fredholm kernel  has attracted increasing attention recently since it  can effectively utilize the data information to improve the prediction
performance. Despite rapid progress on theoretical and experimental evaluations, its generalization analysis has not been explored in learning theory literature. In this paper, we establish the generalization bound of  least square regularized regression with Fredholm kernel, which implies that the fast learning rate $O(l^{-1})$ can be reached under mild capacity conditions. Simulated examples  show that  this Fredholm regression algorithm can  achieve the satisfactory prediction performance.
\end{abstract}


\end{frontmatter}

\section{Introduction}\label{section1}


Inspired from Fredholm integral equations,  Fredholm learning algorithms are designed recently  for density ratio estimation
 \cite{que1} and semi-supervised learning \cite{que2}. Fredholm learning  can be considered as a kernel method with data-dependent  kernel. This kernel usually is called as Fredholm kernel, and can naturally incorporate the data information.    Although its  empirical  performance  has been well demonstrated in the previous works, there is no learning theory analysis on generalization bound and learning rate. It is well known that generalization ability and learning rate are important measures to evaluate the learning algorithm \cite{cucker1,zou1,zou2}. In this paper, we focus on this theoretical theme for regularized least square regression with Fredholm kernel.

In learning theory literature, extensive studies have been established for least square regression with regularized kernel methods, e.g., \cite{shi1,sun1,wu2}.
Although the Fredholm learning in \cite{que2} also can be considered as a regularized kernel method, there are two key features: one is that Fredholm kernel  is associated with the ``inner" kernel and the ``outer" kernel simultaneously, the other is that  for the prediction function is double data-dependent.
These characteristics induce the additional difficulty on learning theory analysis. To overcome the difficulty of generalization analysis,  we introduce novel stepping-stone functions and establish the decomposition on excess generalization error. The generalization bound is estimated in terms of the capacity conditions on the hypothesis spaces associated with the ``inner" kernel and the ``outer" kernel, respectively. In particular, the derived result  implies that fast learning rate with $\mathcal{O}(l^{-1})$
can be reached with proper parameter selection, where $l$ is the number of labeled data. To best of our knowledge, this is the first discussion on generalization error analysis for learning with Fredholm  kernel.

 The rest of this paper is organized as follows. Regression algorithm with Fredholm kernel is introduced in Section \ref{section2} and its generalization analysis is presented in Section \ref{section3}. The  proofs of main results are listed in Section \ref{section4}. Simulated examples are provided in Section \ref{section5} and a brief conclusion is summarized  in Section \ref{section6}.

\section{Regression with Fredholm kernel}\label{section2}

Let $\mathcal{X}\subset\mathbb{R}^d$ be a compact input space and $\mathcal{Y}\subset[-M,M]$ for some constant $M>0$.
The labeled data $\mathbf{z}=\{z_i\}_{i=1}^{l}=\{(x_i,y_i)\}_{i=1}^l$ are drawn independently from a distribution $\rho$ on
$\mathcal{Z}:=\mathcal{X}\times\mathcal{Y}$ and the unlabeled data $\{x_{l+j}\}_{j=1}^u$ are derived random independently according to the marginal distribution $\rho_{\mathcal{X}}$ on $\mathcal{X}$.
Given $\mathbf{z}, \mathbf{x}=\{x_i\}_{i=1}^{l+u}$, the main purpose of semi-supervised regression is to find a good approximation of the regression function
\begin{eqnarray*}
f_{\rho}(x)=\int_{\mathcal{Y}}yd\rho(y|\mathcal{X})=\argmin_f\int_{\mathcal{Z}}(y-f(x))^2d\rho(x,y).
\end{eqnarray*}

In learning theory,
\begin{eqnarray*}
\mathcal{E}(f):=\int_\mathcal{Z}(y-f(x))^2d\rho(x,y)
\end{eqnarray*}
and its discrete version
\begin{eqnarray*}
\mathcal{E}_{\mathbf{z}}(f):=\frac{1}{l}\sum_{i=1}^l(y-f(x_i))^2
\end{eqnarray*}
are called as the expected risk and the empirical risk of function $f:\mathcal{X}\rightarrow\mathbb{R}$, respectively.

Let $w(x,x')$ be a continuous bounded function on $\mathcal{X}^2$ with $\omega:=\sup\limits_{x,x'}w(x,x')<\infty$. Define the integral operator $L_w$ as
\begin{eqnarray*}
L_wf(x)=\int_{\mathcal{X}} w(x,t)f(t)d\rho_{\mathcal{X}}(t),      \forall f\in L_{\rho_{\mathcal{X}}}^2,
\end{eqnarray*}
where $L_{\rho_{\mathcal{X}}}^2$ is the space of square-integrable functions.

Let $\mathcal{H}_K$ be a reproducing kernel Hilbert space (RKHS) associated with Mercer kernel
$K:\mathcal{X}^2\rightarrow\mathbb{R}$. Denote $\|\cdot\|_K$ as the corresponding norm of $\mathcal{H}_K$
and assume the upper bound $\kappa:=\sup\limits_{x,x'\in\mathcal{X}}K(x,x')<\infty$.

If choose
$L_{w}\mathcal{H}=\{L_wf, f\in\mathcal{H}_K\}$ as the hypothesis space, the learning problem can be considered as to solve the Fredhom integral equation $L_wf(x)=y$. Sine the distribution $\rho$ is unknown, we consider the empirical version of $L_wf$
associated with $\mathbf{x}=\{x_i\}_{i=1}^{l+u}$, which is defined as
\begin{eqnarray*}
L_{w,\mathbf{x}}f(x)=\frac{1}{l+u}\sum_{i=1}^{l+u}w(x,x_i)f(x_i).
\end{eqnarray*}

In the Fredholm learning framework, the prediction function is constructed from the data dependent hypothesis space
\begin{eqnarray*}
L_{w,\mathbf{x}}\mathcal{H}=\{L_{w,\mathbf{x}}f, f\in\mathcal{H}_K\}.
\end{eqnarray*}

Given $\mathbf{z}, \mathbf{x}$, least-square regression  with Fredholm kernel (LFK)  can be formulated as the following optimization
\begin{eqnarray}
f_{\mathbf{z}}:=f_{\mathbf{z},\mathbf{x}}=\argmin\limits_{f\in\mathcal{H}_K}\{\mathcal{E}_{\mathbf{z}}(L_{w,\mathbf{x}}f)+\lambda\|f\|_{K}^2\},      \label{algorithm1}
\end{eqnarray}
where $\lambda>0$ is a regularization parameter.

\begin{remark}
Equation (\ref{algorithm1}) can be considered as a discrete and regularized version of the Fredholm integral equation
$L_wf=y$. When $w$ is the $\delta$-function, (\ref{algorithm1}) becomes the regularized least square
regression in RKHS
\begin{eqnarray}
\tilde{f}_{\mathbf{z}}=\argmin\limits_{f\in\mathcal{H}_K}\{\mathcal{E}_{\mathbf{z}}(f)+\lambda\|f\|_{K}^2\}.      \label{algorithm2}
\end{eqnarray}
When $\mathbf{x}=\{x_i\}_{i=1}^l$ and replacing $\|f\|_{K}^2$ with $\sum_{i=1}^l|f(x_i)|^q,q=1,2$, (\ref{algorithm1}) is equivalent to the data-dependent
 coefficient regularization
\begin{eqnarray*}
\tilde{f}_{\mathbf{z}}(x)=\sum_{i=1}^l\alpha_{\mathbf z,i}w(x,x_i),
\end{eqnarray*}
where
\begin{eqnarray}
\alpha_{\mathbf z}
=\argmin\limits_{\alpha\in\mathbb R^l}\Big\{\mathcal{E}_{\mathbf{z}}(\sum_{i=1}^l\alpha_iw(\cdot,x_i))+\lambda\sum_{i=1}^l|
\alpha_i|^q \Big\}. \label{algorithm3}
\end{eqnarray}
It is well known that (\ref{algorithm2}) and (\ref{algorithm3}) have been studied extensively in learning literatures, see, e.g.
\cite{feng3,shi1,sun1}.
These results relied on error analysis techniques for data independent hypothesis space \cite{cucker1,cucker2,zou2} and data dependent hypothesis space \cite{hong3,sun1,sun2,feng3}, respectively. Therefore, the Fredholm learning provides a novel framework for regression related with the data independent space $\mathcal{H}_K$ and the data dependent hypothesis space $L_{w,\mathbf{x}}\mathcal{H}$ simultaneously.
\end{remark}

\begin{remark} Equation (\ref{algorithm1}) involves the ``inner" kernel $K$ and the ``outer" kernel $w$.
Denote $$\hat{K}(x,x')=\frac{1}{(l+u)^2}\sum_{i,j=1}^{l+u}w(x,x_i)K(x_i,x_j)w(x,x_j),$$
$\hat{\mathbf{K}}=(\hat{K}(x_i,x_j))_{i,j=1}^l$,
and $\mathbf{Y}=(y_1,\cdots,y_l)^T$. It has been demonstrated in \cite{que2} that
\begin{eqnarray}
L_{w,\mathbf{x}}f_{\mathbf z}(x)=\frac{1}{l+u}\sum_{i=1}^{l+u}w(x,x_i)f_{\mathbf z}(x_i)=\sum_{s=1}^l\hat{K}(x,x_s)\alpha_s,  \label{algorithm4}
\end{eqnarray}
where $\alpha=(\alpha_1,\cdots,\alpha_l)^T=(\hat{\mathbf{K}}+\lambda I)^{-1}\mathbf{Y}$.
Therefore, Fredholm regression in (\ref{algorithm1}) can be implemented efficiently and the data-dependent kernel $\hat{K}(x,x')$ is called Fredholm kernel in \cite{que2}.
\end{remark}

\section{Generalization bound}\label{section3}

To provide the estimation on the excess risk, we introduce some conditions on the  hypothesis space capacity and the approximation ability of Fredholm learning framework.

For $R>0$, denote
$$
B_R=\{f\in\mathcal{H}_K:\|f\|_K\leq R\}
$$
and
$$
\tilde{B}_R=\{f=\sum\limits_{i=1}^{l+u}\alpha_iw(\cdot,u_i):  \sum\limits_{i=1}^{l+u}|\alpha_i|\leq R,u_i\in\mathcal X\}.$$
  For any $\varepsilon>0$ and function space $\mathcal{F}$, denote $\mathcal{N}_{\infty}(\mathcal{F},   \varepsilon)$ as the covering number with $\ell_{\infty}$-metric.

\begin{assumption}
(Capacity condition) For the ``inner" kernel $K$ and the ``outer" kernel $w$, there exists positive constants $s$ and $p$ such that for any $\varepsilon>0$,
$
\log\mathcal{N}_\infty( B_1,\varepsilon)\leq c_{s,K}\varepsilon^{-s}
$
and
$
\log\mathcal{N}_\infty( \tilde{B}_1,\varepsilon)\leq c_{p,w}\varepsilon^{-p}$,
where $c_{s,K},c_{p,w}>0$ are constants independent of $\varepsilon$.
\label{condition1}
\end{assumption}

It is worthy notice that the capacity condition has been well studied in \cite{cucker1,cucker2,shi1}. In particular, this condition holds true when
setting the ``inner" and ``outer" kernels as Gaussian kernel.

For a function $f:\mathcal X\rightarrow\mathbb R$ and $q\in [1,+\infty)$, denote the $L^q$-norm on $\mathcal X$ as
\begin{eqnarray*}
\|f\|_q:=\|f\|_{L_{\rho_{\mathcal X}}^q}=\Big(\int_{\mathcal X}|f(x)|^qd\rho_{\mathcal{X}}(x)\Big)^{\frac{1}{q}}.
\end{eqnarray*}

Define the data independent regularized function
\begin{eqnarray*}
f_{\lambda}=\argmin\limits_{f\in{\mathcal{H}_K}}\{\|L_{w}f-f_{\rho}\|_2^2+\lambda\|f\|_K^2\}.
\end{eqnarray*}

The predictor associated with $f_\lambda$ is
\begin{eqnarray*}
L_wf_\lambda=\int_{\mathcal{X}}w(x,t)f_{\lambda}(t)d\rho_{\mathcal X}(t)
\end{eqnarray*}
and the approximation ability of Fredholm scheme in $\mathcal{H}_K$ is characterized by
\begin{eqnarray*}
D(\lambda)=\mathcal E(L_wf_{\lambda})-\mathcal E(f_{\rho})+\lambda\|f_{\lambda}\|_K^2.
\end{eqnarray*}

\begin{assumption}(Approximation condition) There exists a constant $\beta\in(0,1]$ such that
\begin{eqnarray*}
D(\lambda)\leq c_\beta\lambda^{\beta},~~ \forall\lambda>0,
\end{eqnarray*}
where $c_\beta$ is a positive constant independent of $\lambda$.\label{condition2}
\end{assumption}

This approximation condition relies on the regularity of $f_\rho$, and has been investigated extensively in \cite{cucker2,sun1,hong4}.
To get tight estimation, we introduce the projection operator
\begin{eqnarray*}
\pi(f)(x)=
\left\{
                       \begin{array}{ll}
                        M, & \hbox{if $f(x)>M$;} \\
                         f(x), & \hbox{if $|f(x)|\leq M$;} \\
                         -M, & \hbox{if $f(x)<-M$.}
                       \end{array}
                     \right.
\end{eqnarray*}

It is a position to present the generalization bound.

\begin{theorem} Under Assumptions \ref{condition1} and \ref{condition2} , there exists
\begin{eqnarray*}
\mathcal E(\pi(L_{w,\mathbf{x}}f_{\mathbf z}))-\mathcal E(f_\rho)
\leq c\log^2(6/{\delta})(\lambda^{-\frac{s}{2+s}}l^{-\frac{s}{2+s}}+\lambda^{\beta}
+\lambda^{\beta-1}l^{-\frac{s}{2+p}}),
\end{eqnarray*}
where $c$ is a positive constant independent of $l,\lambda,\delta$
 \label{theorem1}
\end{theorem}

The generalization bound  in Theorem \ref{theorem1} depends on the capacity condition , the approximation condition, the regularization parameter $\lambda$, and the number of labeled data. In particular, the labeled data is the key factor on the excess risk without the additional assumption on the marginal distribution. This observation is consistent with the previous  analysis for semi-supervised learning  \cite{belkin,hong1}.

To understand the learning rate of Fredholm regression, we present the following result where $\lambda$ is chosen properly.

\begin{theorem} Under Assumptions \ref{condition1} and \ref{condition2}, for any $0<\delta<1$, with
confidence $1-\delta$, there exists some positive constant $\tilde c$ such that
\begin{eqnarray*}
\mathcal E(\pi(L_{w,\mathbf{x}}f_{\mathbf z}))-\mathcal E(f_\rho)\leq \tilde c\log^2(6/\delta)l^{-\theta},
\end{eqnarray*}
where
\begin{eqnarray*}
\theta=
\left\{
                       \begin{array}{ll}
                        \min\{\frac{2\beta}{2+p},\frac{2}{2+s}-\frac{2s}{(2+s)(2+p)}\}, & \hbox{ $\lambda=l^{-\frac{2}{2+p}}$;} \\
                        \min\{\frac{2\beta}{2\beta+s\beta+s},\frac{(2\beta+s\beta+s)(\beta-1)}{2+s}-\frac{2}{2+p}\}, & \hbox{ $\lambda=l^{-\frac{2}{2\beta+s\beta+s}}$.} \\
                       \end{array}
                     \right.
\end{eqnarray*}
\label{theorem2}
\end{theorem}

Theorem \ref{theorem2} tells us that Fredholm regression has the learning rate with polynomial decay. When $s=p$, there exists some constant $\bar{c}>0$ such that
\begin{eqnarray*}
\mathcal E(\pi(L_{w,\mathbf{x}}f_{\mathbf z}))-\mathcal E(f_\rho)\leq \bar c\log(6/\delta)l^{-\theta}
\end{eqnarray*}
with confidence $1-\delta$,
where
\begin{eqnarray*}
\theta=
\left\{
                      \begin{array}{ll}
                       \frac{2\beta}{2+s},& \beta\in(0,\frac{2}{2+s}]; \\
                       \frac{2\beta}{s+2\beta+s\beta},& \beta\in(\frac{2}{2+s},+\infty]. \\
                       \end{array}
                     \right.
\end{eqnarray*}
and the rate is derived by setting
\begin{eqnarray*}
\lambda=
\left\{
                      \begin{array}{ll}
                       l^{-\frac{2}{2+s}},& \beta\in(0,\frac{2}{2+s}]; \\
                       l^{-\frac{2}{s+2\beta+s\beta}},& \beta\in(\frac{2}{2+s},+\infty]. \\
                       \end{array}
                     \right.
\end{eqnarray*}
This  learning rate can be arbitrarily close to $\mathcal O(l^{-1})$ as $s$ tends to zero,
which is regarded as the fastest learning rate for regularized regression in the learning theory literature. This result verifies the LFK in  (\ref{algorithm1}) inherits the theoretical characteristics of least square regularized
regression in RKHS \cite{cucker2,wu2} and in data dependent hypothesis spaces \cite{shi1,sun2}.

\section{Error analysis}\label{section4}

We first present the decomposition on the excess risk $\mathcal E(\pi(L_{w,\mathbf{x}}f_{\mathbf z}))-\mathcal E(f_\rho)$, and then establish the upper bounds of different error terms.

\subsection{Error decomposition}\label{section4.1}
According to the definitions of $f_{\mathbf z}, f_\lambda$, we can get the following error decomposition.

\begin{proposition}\label{proposition1}
For $f_\mathbf z$ defined in (\ref{algorithm1}), there holds
\begin{eqnarray*}
\mathcal E(\pi(L_{w,\mathbf{x}}f_{\mathbf z}))-\mathcal E(f_\rho)\leq E_1+E_2+E_3+D(\lambda),
\end{eqnarray*}
where
\begin{eqnarray*}
E_1&=&\mathcal E(\pi(L_{w,\mathbf{x}}f_{\mathbf z}))-\mathcal E(f_\rho)-(\mathcal E_{\mathbf z}(\pi(L_{w,\mathbf{x}}f_{\mathbf z}))-\mathcal E_{\mathbf z}(f_\rho)),\\
E_2&=&\mathcal E_{\mathbf z}(L_{w,\mathbf{x}}f_{\lambda})-\mathcal E_{\mathbf z}(f_\rho)-(\mathcal E(L_{w,\mathbf{x}}f_{\lambda})-\mathcal E(f_\rho)),
\end{eqnarray*}
and
$$
E_3=\mathcal E(L_{w,\mathbf{x}}f_{\lambda})-\mathcal E(L_{w}f_{\lambda}).$$
\end{proposition}
{\bf {Proof}}:
By introducing the middle function $L_{w,\mathbf{x}}f_{\lambda}$, we get
\begin{eqnarray*}
&&\mathcal E(\pi(L_{w,\mathbf{x}}f_{\mathbf z}))-\mathcal E(f_\rho)\\
&\leq& \mathcal E(\pi(L_{w,\mathbf{x}}f_{\mathbf z}))-\mathcal E_{\mathbf z}(\pi(L_{w,\mathbf{x}}f_{\mathbf z}))
+[\mathcal E_{\mathbf z}(L_{w,\mathbf{x}}f_{\mathbf z})+\lambda\|f_{\mathbf z}\|_K^2-(\mathcal E_{\mathbf z}(L_{w,\mathbf{x}}f_{\lambda})+\lambda\|f_\lambda\|_K^2)]\\
&&+\mathcal E_{\mathbf z}(L_{w,\mathbf{x}}f_{\lambda})-\mathcal E(L_{w,\mathbf{x}}f_{\lambda})+\mathcal E(L_{w,\mathbf{x}}f_{\lambda})-\mathcal E(L_{w}f_{\lambda})
+\mathcal E(L_{w}f_{\lambda})-\mathcal E(f_\rho)+\lambda\|f_{\lambda}\|_K^2 \\
&\leq& E_1+E_2+E_3+D(\lambda)
\end{eqnarray*}
where the last inequality follows from the definition $f_{\mathbf z}$.
This completes the proof.$\blacksquare$

In learning theory, $E_1,E_2$ are  called the sample error, which describe the difference between the empirical risk and the expected risk. $E_3$ is  called the hypothesis error which reflects the divergence  of expected risks between the data independent function $L_wf_\lambda$ and data dependent function $L_{w,\mathbf{x}}f_\lambda$.

\subsection{Estimates of sample error}\label{subsection4.2}
We introduce the concentration inequality in \cite{wu1} to measure the divergence between the empirical risk and the expected risk.

\begin{lemma}Let $\mathcal F$ be a measurable function set on $\mathcal Z$. Assume that,   for any $f\in\mathcal F$,   $\|f\|_{\infty}\leq B$ and $E(f^2)\leq cEf$ for some positive constants $B,   c$. If for some $a>0$
and $s\in(0,   2)$,   $\log \mathcal{N}_2(\mathcal F,   \varepsilon)\leq a\varepsilon^{-s}$
for any $\varepsilon>0$,   then there exists a constant $c_s$ such that for any
$\delta\in(0,   1)$,
\begin{eqnarray*}
\Big|Ef-\frac{1}{m}\sum_{i=1}^mf(z_i)\Big|
\leq c_s\max\{c^{\frac{2-s}{2+s}},   B^{\frac{2-s}{2+s}}\}
(\frac{a}{m})^{\frac{2}{2+s}}
+\frac{1}{2}Ef+\frac{(2c+18B)\log(1/\delta)}{m}
\end{eqnarray*}
with confidence at least $1-2\delta$.
\label{lemma1}
\end{lemma}

To estimate $E_1$, we consider the function set containing $f_{\mathbf z}$ for any $\mathbf{z}\in \mathcal{Z}^l$, $\mathbf{u}\in \mathcal{X}^u$.
The definition $f_{\mathbf z}$ in (\ref{algorithm1}) tells us that $\|f_{\mathbf z}\|_K\leq\frac{M}{\sqrt\lambda}$.
Hence, $\forall \mathbf{z}\in \mathcal{Z}^l,f_{\mathbf z}\in B_R$ with $R=\frac{M}{\sqrt\lambda}$ and $\|f_{\mathbf z}\|_{\infty}\leq\frac{\kappa M}{\sqrt\lambda}$.

\begin{proposition}\label{proposition2}
Under Assumption \ref{condition1}, for any $0<\delta<1$,
\begin{eqnarray*}
E_1 \leq\frac{1}{2}(\mathcal E(\pi(L_{w,\mathbf{x}}f_{\mathbf z}))-\mathcal E(f_\rho))+c_1\lambda^{-\frac{s}{2+s}}m^{-\frac{s}{2+s}}
+176M^2l^{-1}\log(1/\delta)
\end{eqnarray*}
with confidence $1-\delta$.
\end{proposition}
{\bf{Proof}}:
For $f\in B_R,\mathbf{z}\in \mathcal{Z}^l, \mathbf{x}\in\mathcal \mathcal{X}^{l+u}$, denote
\begin{eqnarray*}
G_R=\{g(z)=(y-\pi(L_{w,\mathbf{x}}f))^2-(y-f_\rho(x))^2\}.
\end{eqnarray*}
For any $z\in \mathcal{Z}$,
\begin{eqnarray*}
|g(z)|\leq|2y-\pi(L_{w,\mathbf{x}}f)(x)-f_\rho(x)||\pi(L_{w,\mathbf{x}}f)-f_\rho(x)|
\leq 8M^2.
\end{eqnarray*}
Moreover,
\begin{eqnarray*}
Eg^2\leq 16M^2E(\pi(L_{w,\mathbf{x}}f)(x)-f_\rho(x))^2=16M^2Eg.
\end{eqnarray*}
For any $f_1,f_2\in B_R$,  there exists
\begin{eqnarray*}
|g_1(z)-g_2(z)|
\leq \frac{4M}{l+u}\Big|\sum_{i=1}^{l+u}(f_1(x_i)-f_2(x_i))w(x,x_i)\Big|
\leq 4M\omega\|f_1-f_2\|_{\infty}.
\end{eqnarray*}
This relation implies that
\begin{eqnarray*}
\log\mathcal{N}_{\infty}( G_R,\varepsilon)\leq\log\mathcal{N}_{\infty}(B_1,\frac{\varepsilon}{4M\omega R})\leq c_{s,K}(4M\omega R)^s\varepsilon^{-s},
\end{eqnarray*}
where the last inequality from Assumption \ref{condition1}.

Applying the above estimates  to Lemma  \ref{lemma1}, we derive that
\begin{eqnarray*}
&& Eg-\frac{1}{l}\sum_{i=1}^{l}g(z_i)\\
&\leq& \frac{1}{2}Eg+\max\{16M^2\omega,8M^2\}^{\frac{2-s}{2+s}}
 c_{s,K}^{\frac{2}{2+s}}(4M\omega R)^{\frac{2s}{2+s}}l^{-\frac{2}{2+s}}
+ 176M^2l^{-1}\log(1/\delta)
\end{eqnarray*}
with confidence $1-\delta$.

Considering $f_\mathbf{z}\in B_R$ with $R=\frac{M}{\sqrt\lambda}$, we obtain the desired result.
$\blacksquare$

\begin{proposition}\label{proposition3}
Under Assumption 1, with confidence $1-4\delta$, there holds
\begin{eqnarray*}
E_2\leq\frac{1}{2}E_3+\frac{1}{2}D(\lambda)+c_2D(\lambda)\lambda^{-1}l^{-\frac{2}{2+p}}
\log(1/\delta),
\end{eqnarray*}
where $c_2$ is a positive constant independent of $\lambda, m, \delta$.
\end{proposition}
{\bf{Proof}}:
Denote
\begin{eqnarray*}
\mathcal{G}=\{g_{\mathbf{v},\lambda}:g_{\mathbf{v},\lambda}(x)=L_{w,\mathbf{v}}f_{\lambda}(x), x,v_i\in\mathcal{X}\}.
\end{eqnarray*}

From the definition $f_{\lambda}$, we can deduce that $\forall g\in\mathcal{G}, g\in\tilde{B}_R$ with $R=\omega\kappa\sqrt{\frac{D(\lambda)}{\lambda}}$. For $z\in \mathcal{Z}, \mathbf{v}\in\mathcal{X}^{l+u}$, define
$$
\mathcal{H}=\{h(z)=(y-L_{w,\mathbf{v}}f_{\lambda}(x))^2-(y-f_\rho(x))^2\}.$$
It is easy to check that for any $z\in \mathcal{Z}$
\begin{eqnarray}
|h(z)|&=&|2y-L_{w,\mathbf{v}}f_{\lambda}(x)-f_\rho(x)|\cdot|L_{w,\mathbf{v}}f_{\lambda}(x)-f_\rho(x)| \nonumber\\
&\leq& (3M+\omega\|f_\lambda\|_{\infty})^2\leq\Big(3M+\omega\kappa\sqrt{\frac{D(\lambda)}{\lambda}}\Big)^2.  \label{p11}
\end{eqnarray}

Then,
\begin{eqnarray}
Eh^2&=&E(2y-L_{w,\mathbf{v}}f_{\lambda}(x)-f_\rho(x))^2(L_{w,\mathbf{v}}f_{\lambda}(x)-f_\rho(x))^2 \nonumber\\
&\leq&\Big(3M+wk\sqrt{\frac{D(\lambda)}{\lambda}}\Big)^2Eh. \label{p22}
\end{eqnarray}

For any $\mathbf{u},\mathbf{v}\in\mathcal{X}^{l+u}$, there exists
\begin{eqnarray*}
\|h_1-h_2\|_{\infty}
&=&\sup_{z}|(y-L_{w,\mathbf{u}}f_{\lambda}(x))^2-(y-L_{w,\mathbf{v}}f_{\lambda}(x))^2| \\
&\leq&2\Big(M+\omega\kappa\sqrt{\frac{D(\lambda)}{\lambda}}\Big)\|L_{w,\mathbf{u}}f_{\lambda}-L_{w,\mathbf{v}}f_{\lambda}\|_{\infty} \\
&=&2\Big(M+\omega\kappa\sqrt{\frac{D(\lambda)}{\lambda}}\Big)\|g_{\mathbf{u},\lambda}-g_{\mathbf{v},\lambda}\|_{\infty}.
\end{eqnarray*}

Then from Assumption \ref{condition1},
\begin{eqnarray}
\log\mathcal{N}_{\infty}(\mathcal{H},\varepsilon)
\leq\log\mathcal{N}_{\infty}
\Big(\tilde{B}_R,\frac{\varepsilon}{2(M+\omega\kappa\sqrt{\frac{D(\lambda)}{\lambda}})}\Big)
\leq4c_{p,w}\Big(M+\omega\kappa\sqrt{\frac{D(\lambda)}{\lambda}}\Big)^{2p}
\varepsilon^{-p}. \label{p33}
\end{eqnarray}

Combining  (\ref{p11})-(\ref{p33}) with Lemma  \ref{lemma1}, we get with confidence $1-\delta$
\begin{eqnarray*}
 E_2\leq\frac{1}{2}(\mathcal E(L_{w,\mathbf{x}}f_{\lambda})-f_\rho)+ (M+\omega\kappa\sqrt{\frac{D(\lambda)}{\lambda}})^2l^{-\frac{2}{2+p}}\\
\cdot c_p(4c_{p,w})^{\frac{2}{2+p}}+\frac{20(3M+\omega\kappa\sqrt{\frac{D(\lambda)}{\lambda}})\log(1/\delta)}{l}.
\end{eqnarray*}

Considering $
 \mathcal E(L_{w,\mathbf{x}}f_{\lambda})-\mathcal E(f_\rho)
 \leq E_3+D(\lambda)$, we get the desired result.$\blacksquare$

\subsection{Estimate of hypothesis error}\label{section4.3}
The following  concentration inequality with values in Hilbert space can be found in \cite{pinelis}, which is used in our analysis.

\begin{lemma}\label{lemma2}
Let $\mathcal H$ be a Hilbert space and $\xi$ be independent random variable on $Z$ with
values in $\mathcal H$. Assume that $\|\xi\|_{\mathcal H}\leq \tilde M<\infty$ almost surely.
Let $\{z_i\}_{i=1}^m$ be independent random samples from $\rho$.
Then, for any $\delta\in(0,1)$,
\begin{eqnarray*}
 \Big\|\frac{1}{m}\sum_{i=1}^m\xi(z_i)-E\xi\Big\|_{\mathcal H}\leq\frac{2\tilde M\log(\frac{1}{\delta})}{M}
 +\sqrt{\frac{2E\|\xi\|_{\mathcal H}^2\log(\frac{1}{\delta})}{m}}
\end{eqnarray*}
holds true with confidence $1-\delta$.
\end{lemma}

Now we turn to estimate $E_3$, which reflects the affect of inputs $\mathbf{x}=\{x_i\}_{i=1}^{l+u}$
to the regularization function $f_{\lambda}$.

\begin{proposition}\label{proposition}
For any $0<\delta<1$, with confidence $1-\delta$, there holds
\begin{eqnarray*}
E_3\leq 24\omega^2\kappa^2\log^2(\frac{1}{\delta})D(\lambda)\lambda^{-1}(l+u)^{-1}+D(\lambda).
\end{eqnarray*}
\end{proposition}
{\bf{Proof}}:
Note that
\begin{eqnarray}
&&\mathcal E(L_{w,\mathbf{x}}f_{\lambda})-\mathcal E(L_{w}f_{\lambda})
\nonumber \\
&\leq& \|L_{w,\mathbf{x}}f_{\lambda}-L_{w}f_{\lambda}\|_2\cdot(\|L_{w,\mathbf{x}}f_{\lambda}-f_{\rho}\|_2+
\|L_{w}f_{\lambda}-f_{\rho}\|_2)\nonumber\\
&\leq& \|L_{w,\mathbf{x}}f_{\lambda}-L_{w}f_{\lambda}\|_2(\|L_{w,\mathbf{x}}f_{\lambda}-L_{w}f_{\lambda}\|_2+
2\|L_{w}f_{\lambda}-f_{\rho}\|_2)\nonumber\\
&\leq& 2\|L_{w,\mathbf{x}}f_{\lambda}-L_{w}f_{\lambda}\|_2^2+\|L_{w}f_{\lambda}-f_{\rho}\|_2^2\nonumber\\
&\leq& 2\|L_{w,\mathbf{x}}f_{\lambda}-L_{w}f_{\lambda}\|_2^2+D(\lambda). \label{p111}
\end{eqnarray}

Denote $\xi(x_i)=f_{\lambda}(x_i)w(\cdot,x_i)$, which is continuous and bounded function on
$\mathcal X$. Then
\begin{eqnarray*}
L_{w,\mathbf{x}}f_{\lambda}=\frac{1}{l+u}\sum_{i=1}^{l+u}\xi(x_i)
\end{eqnarray*}
and
\begin{eqnarray*}
L_{w}f_{\lambda}=\int w(\cdot,t)f_{\lambda}(t)d\rho_{\mathcal X}(t)=E\xi.
\end{eqnarray*}

We can deduce that $\|\xi\|_2\leq \omega\|f_{\lambda}\|_{\infty}\leq \omega\kappa\|f_{\lambda}\|_K$ and
$E\|\xi\|_2^2\leq \omega^2\kappa^2\|f_{\lambda}\|_K^2$.
From Lemma \ref{lemma2}, for any $0<\delta<1$, there holds with confidence $1-\delta$
\begin{eqnarray}
\|L_{w,\mathbf{x}}f_{\lambda}-L_{w}f_{\lambda}\|_2
\leq\frac{2\omega\kappa\|f_{\lambda}\|_K\log(\frac{1}{\delta})}{l+u}
+\sqrt{\frac{2\log(\frac{1}{\delta})}{l+u}}\omega\kappa\|f_{\lambda}\|_K. \label{p222}
\end{eqnarray}

Combining (\ref{p111}) and (\ref{p222}), we get with confidence $1-\delta$,
\begin{eqnarray*}
E_3&\leq& 2(\frac{2\omega\kappa\|f_\lambda\|_K\log(\frac{1}{\delta})}{l+u}+\omega\kappa\|f_\lambda\|_K\sqrt
{\frac{2\log(\frac{1}{\delta})}{l+u}})^2+D(\lambda)\\
&\leq&\frac{16\omega^2\kappa^2\|f_\lambda\|_K^2\log^2(\frac{1}{\delta})}{(l+u)^2}
+\frac{8\omega^2\kappa^2\|f_\lambda\|_K^2\log(\frac{1}{\delta})}{l+u}+D(\lambda).
\end{eqnarray*}

Then, the desired result follows from $\|f_\lambda\|_K^2\leq\frac{D(\lambda)}{\lambda}$.
$\blacksquare$

\subsection{Proofs of Theorem 1 and 2}\label{section4.4}

{\bf Proof of Theorem 1:} Combining the estimations in Propositions 1-4, we get with confidence $1-6\delta$,
\begin{eqnarray*}
&&\mathcal E(\pi(L_{w,\mathbf{x}}f_\mathbf z))-\mathcal E(f_\rho)\\
&\leq&\frac{1}{2}(\mathcal E(\pi(L_{w,\mathbf{x}}f_\mathbf z))-\mathcal E(f_\rho))+c_1\lambda^{-\frac{s}{2+s}}l^{-\frac{2}{2+s}}
+176M^2l^{-1}\log(\frac{1}{\delta})\\
&&+ 3D(\lambda)+c_2D(\lambda)\lambda^{-1}l^{-\frac{2}{2+p}}\log(\frac{1}{\delta})
+\frac{36w^2k^2\log^2(\frac{1}{\delta})}{l+u}\frac{D(\lambda)}{\lambda}.
\end{eqnarray*}

Considering $u>0$, for $0<\delta<1$, we have with confidence $1-6\delta$
\begin{eqnarray*}
\mathcal E(\pi(L_{w,\mathbf{x}}f_\mathbf z))-\mathcal E(f_\rho)
\leq c\log^2(\frac{1}{\delta})
[\lambda^{-\frac{s}{2+s}}l^{-\frac{s}{2+s}}+\lambda^{\beta}+\lambda^{\beta-1}l^{-\frac{2}{2+p}}],
\end{eqnarray*}
where $c$ is a constant independent of $l,\lambda,\delta$.

{\bf Proof of Theorem 2:} When setting $\lambda^{\beta}=\lambda^{\beta-1}l^{-\frac{2}{2+p}}$, we obtain
$\lambda=l^{-\frac{2}{2+p}}$. Then, Theorem \ref{theorem1} implies that
\begin{eqnarray*}
\mathcal E(\pi(L_{w,\mathbf{x}}f_\mathbf z))-\mathcal E(f_\rho)
\leq 3c\log^2(\frac{1}{\delta})
l^{-\min\{\frac{2\beta}{2+p},\frac{2}{2+s}-\frac{2s}{(2+s)(2+p)}\}}.
\end{eqnarray*}

When setting $\lambda^{\beta}=\lambda^{-\frac{s}{2+s}}l^{-\frac{2}{2+s}}$, we get
$\lambda=l^{-\frac{2}{2\beta+s\beta+s}}$. Then, with confidence $1-6\delta$
\begin{eqnarray*}
\mathcal E(\pi(L_{w,\mathbf{x}}f_\mathbf z))-\mathcal E(f_\rho)
\leq 3c\log^2(\frac{1}{\delta})
l^{-\min\{\frac{2\beta}{2\beta+s\beta+s},\frac{(2\beta+s\beta+s)(\beta-1)}{2+s}
+\frac{2}{2+p}\}}.
\end{eqnarray*}
This complete the proof of Theorem 2.

\section{Empirical studies }\label{section5}
To verify the effectiveness of LFK in (\ref{algorithm1}), we present some  simulated examples for the  regression problem.  The competing method is support vector machine regression (SVM), which has been used extensively used in  machine learning community (https://www.csie.ntu.edu.tw/~cjlin/libsvm/).
The Gaussian kernel $K(x,t)=\exp\{-\frac{\|x-t\|_2^2}{2\sigma^2}\}$ is used for SVM. For LFK in (\ref{algorithm1}), we consider the following ``inner''  and ``outer'' kernels:
\begin{itemize}
 \item LFK1: $w(x,z)=x^Tz$ and $K(x,z)=\exp\{-\frac{\|x-t\|_2^2}{\sigma^2}\}$.

 \item LFK2: $w(x,z)=\exp\{-\frac{\|x-t\|_2^2}{\sigma^2}\}$ and $K(x,z)=x^Tz$.

  \item LFK3: $w(x,z)=\exp\{-\frac{\|x-t\|_2^2}{\sigma^2}\}$ and  $K(x,z)=\exp\{-\frac{\|x-t\|_2^2}{\sigma^2}\}$.
\end{itemize}
Here the scale parameter $\sigma$ belongs to
$[2^{-5}:2:2^5]$ and the regularization parameter belongs to $[10^{-5}:10:10^5]$ for LFK and SVM.
These parameters are selected  by 4-fold cross validation in this section.

The following  functions are used to generate the simulated data:
\begin{eqnarray*}
f_1(x)&=&sin\Big(\frac{9\pi}{0.35x+1}\Big),~~x\in[0,10]\\
f_2(x)&=&xcos(x),~~x\in[0,10]\\
f_3(x)&=&\min(2|x|-1,1),~~x\in[-2,2]\\
f_4(x)&=&sign(x),~~x\in[-3,3].\\
\end{eqnarray*}
Note that $f_1$ is highly oscillatory,  $f_2$ is smooth, $f_3$ is continuous not smooth, and
$f_4$ is not even continuous. These functions have been used to evaluate regression algorithms in \cite{sun2}.
\begin{table}
\caption{MSE$\pm$STD for LFK and SVM with 50 and 100 training samples}
\begin{center}
\begin{tabular}
{c|ccccc}
\hline
Function     &Number   &SVM   &LFK1     &LFK2  &LFK3
\\
\hline
$f_1$               &50	&$0.041\pm0.033$ 	&$0.434\pm0.032$ 	&$0.423\pm0.059$ 	&$\mathbf{0.036\pm0.053}$
        \\
        &300	&$0.044\pm0.006$ 	&$0.419\pm0.023$ 	&$0.404\pm0.021$ 	&$\mathbf{0.042\pm0.006}$
                             \\
\cline{1-6}
$f_2$              &50	&$0.075\pm0.046$ 	&$18.52\pm1.30$ 	&$18.7\pm1.30$ 	&$\mathbf{0.060\pm0.028}$
        \\
                    &300	&$\mathbf{0.011}\pm0.006$ 	&$17.10\pm0.941$ 	&$17.00\pm1.35$ 	&$0.012\pm\mathbf{0.004}$
         \\
\cline{1-6}					
$f_3$              &50	&$0.013\pm 0.012$ 	&$0.670\pm0.034$ 	&$0.458\pm0.082$ 	&$\mathbf{0.010\pm0.005}$
         \\
                    &300	&$0.004\pm0.001$ 	&$0.667\pm0.013$ 	&$0.427\pm0.020$ 	&$\mathbf{0.003\pm0.001}$
         \\
\cline{1-6}
$f_4$               &50	&$0.076\pm0.027$ 	&$0.260\pm\mathbf{0.012}$ 	&$0.194\pm0.040$ 	&$\mathbf{0.073}\pm0.021$
        \\
                    &300	&$0.039\pm0.018$ 	&$0.251\pm0.017$ 	&$0.158\pm0.026$ 	&$\mathbf{0.032\pm0.009}$
        \\
    \hline
\end{tabular}
\end{center}
\label{tab1}
\end{table}

\begin{figure}[ht]
  \centering
   \subfigure[]{
    \nonumber
   \includegraphics[width=5.8cm]{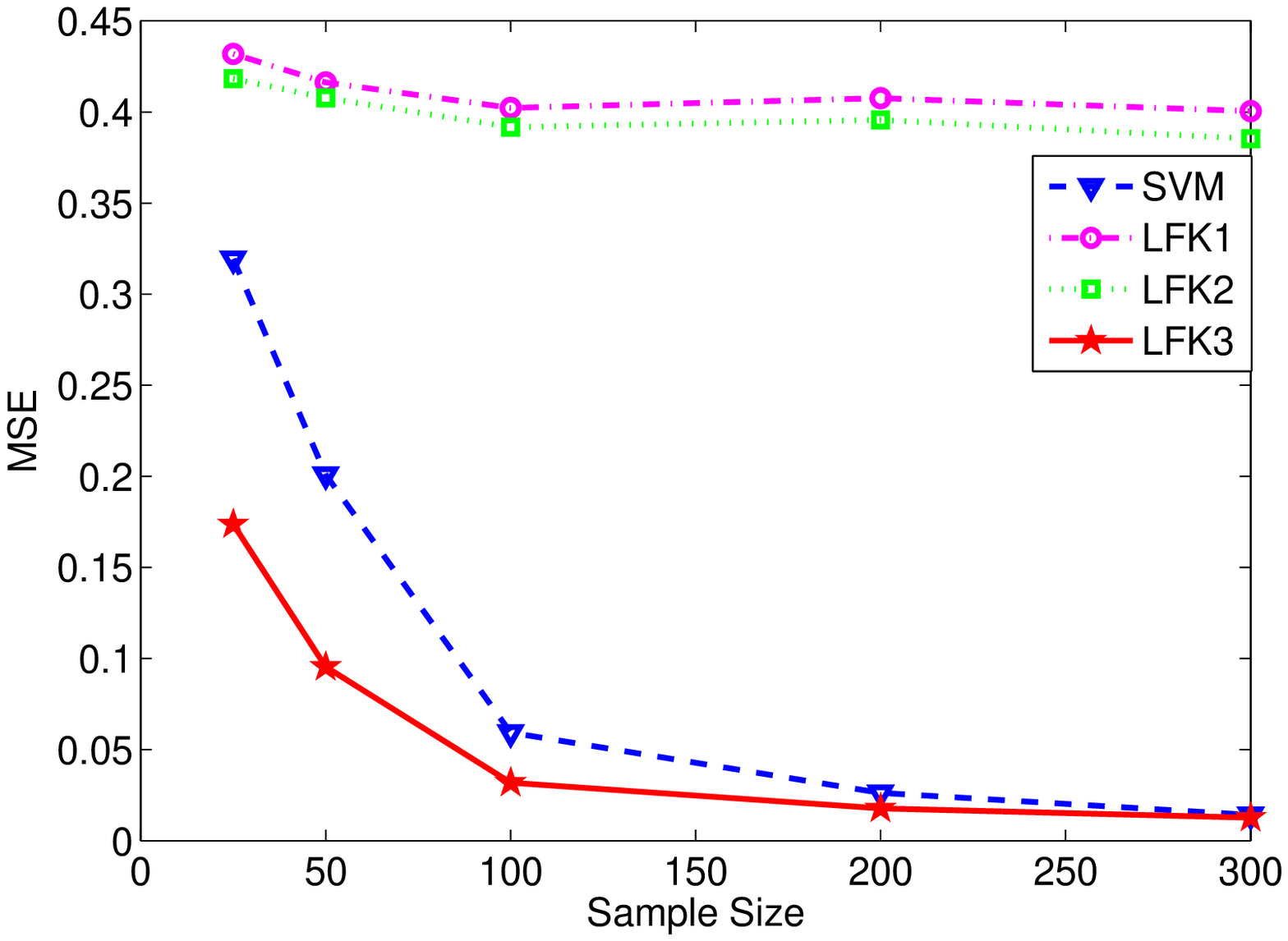}}
      \subfigure[]{
    \nonumber
   \includegraphics[width=5.8cm]{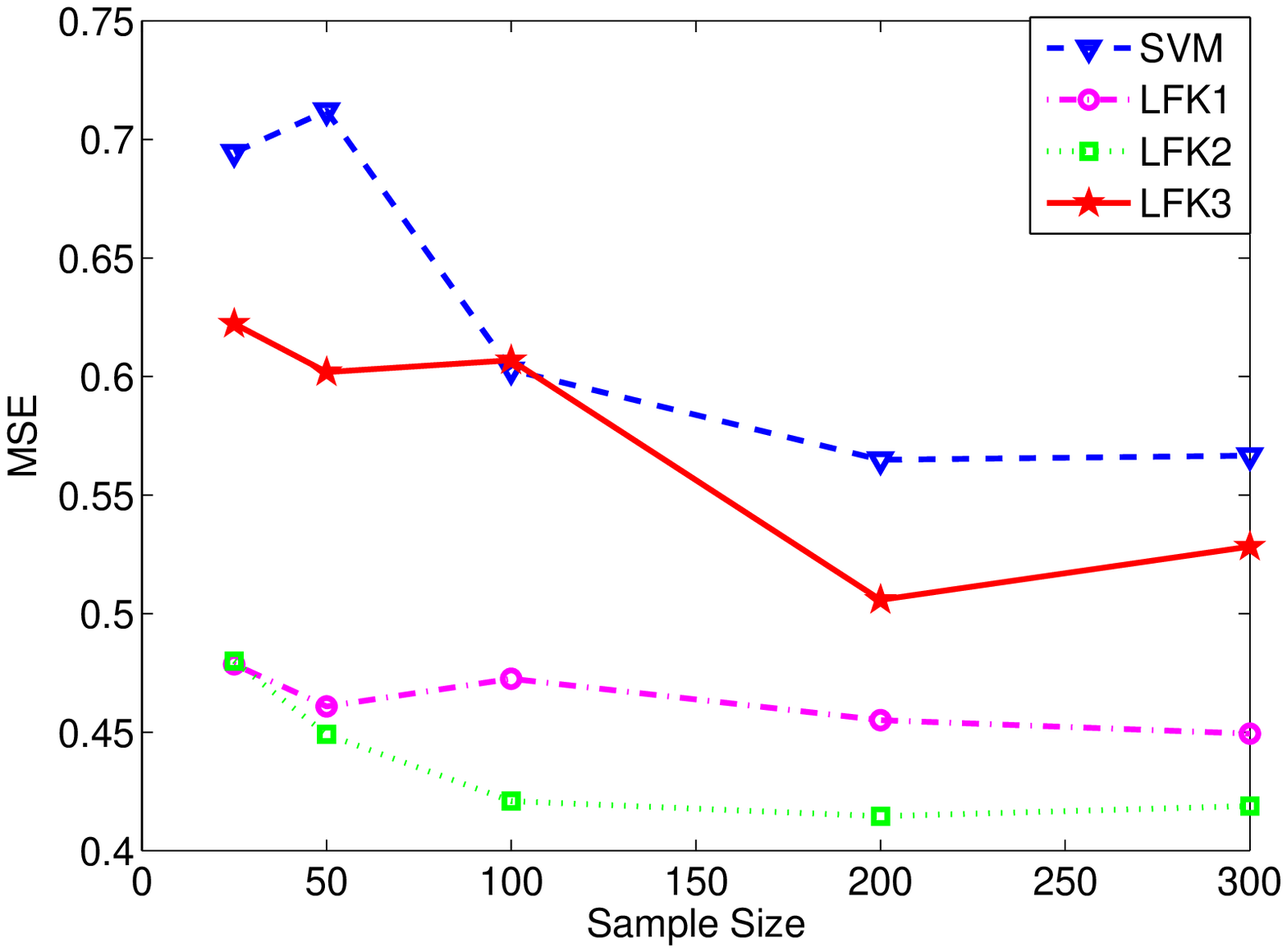}}
      \subfigure[]{
    \nonumber
   \includegraphics[width=5.8cm]{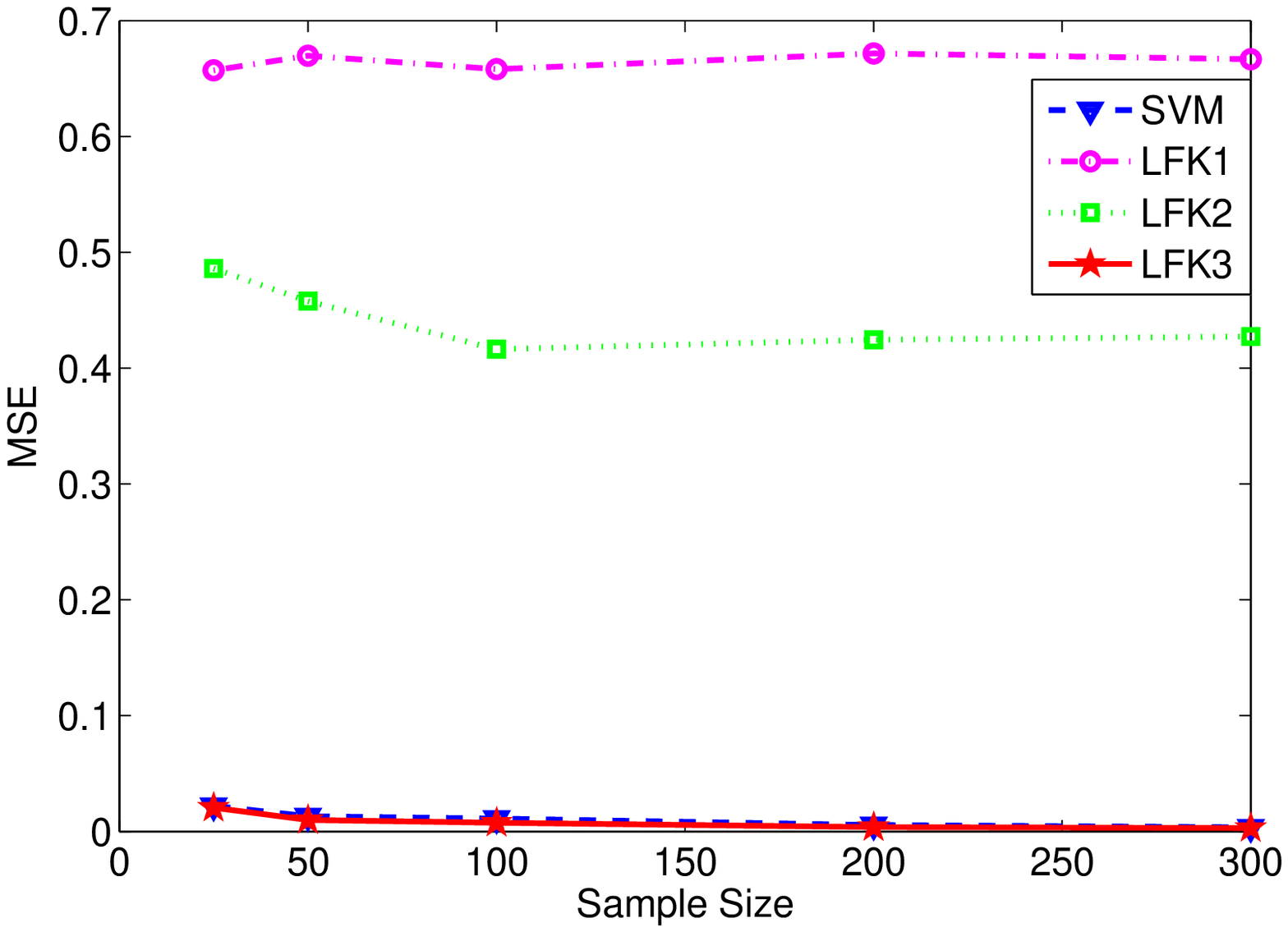}}
        \subfigure[]{
    \nonumber
   \includegraphics[width=5.8cm]{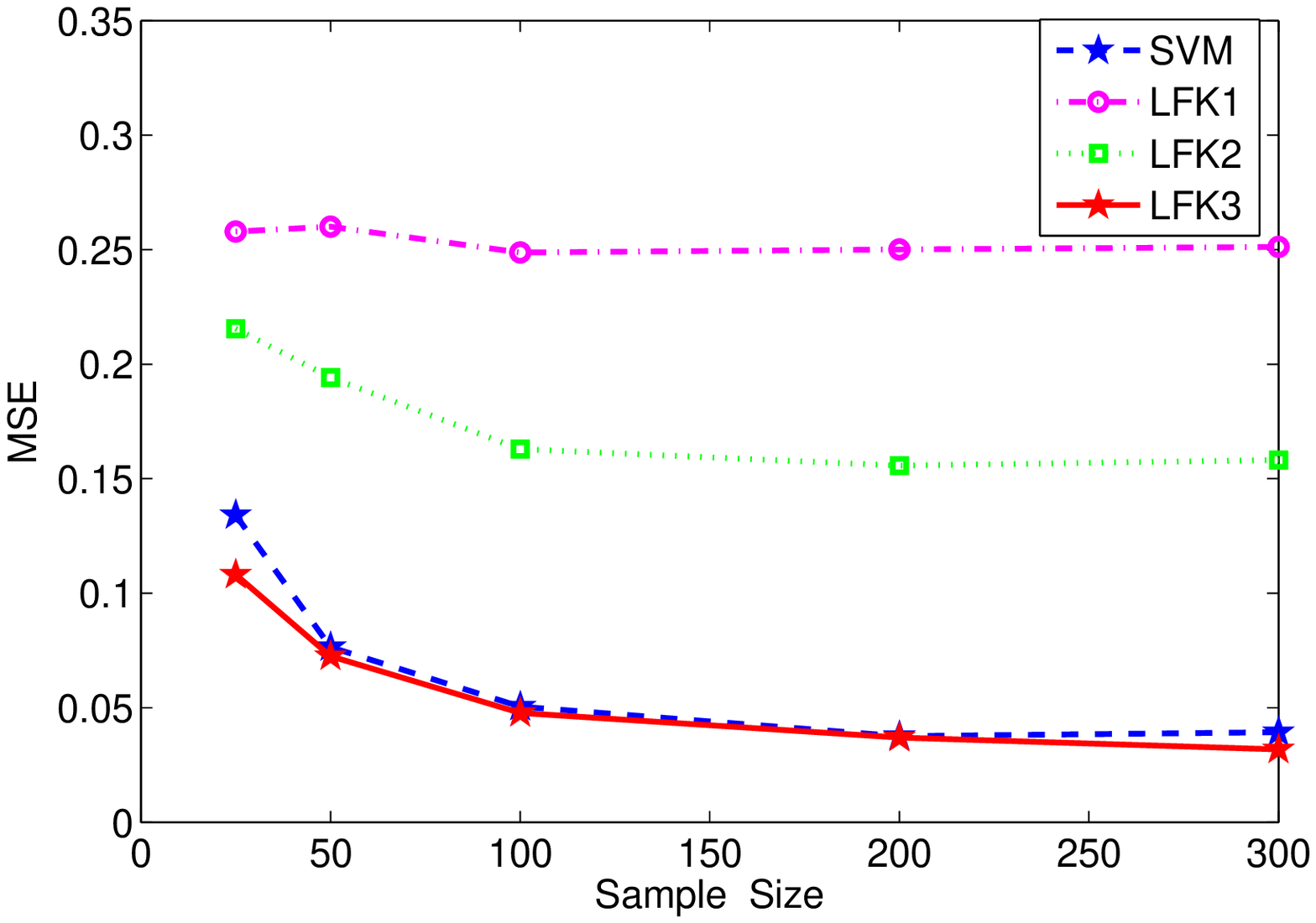}}
    \caption{
  MSE for  Gaussian noise and varying training sample size.
 (a) $f_1$ ; (b) $f_2$; (c) $f_3$; (d) $f_4$.}
\label{RMSE-number1}
\end{figure}

In our experiment, Gaussian noise $N(0,0.01)$ is added to the data respectively.
In each test,  we first draw randomly 1000 samples according to the function and noise distribution, and then obtain  a training set randomly  with sizes  $25, 50, 100, 200, 300$ respectively. Three hundred samples are selected randomly as the test set.
The \emph{Mean Squared Error} (MSE)  is used to evaluate the regression results on synthetic data.
To make the results more convincing, each test is repeated 10 times.  Table \ref{tab1} reports the average MSE and \emph{Standard Deviation} (STD) with 50 training samples and 300 training samples respectively.
Furthermore, we study the impact of the number of training samples  on
the final regression performance. Figure 1  shows the MSE for learning $f_1-f_4$
with  numbers of training samples. These results illustrate that LFK   has competitive performance compared with SVM.

\section{Conclusion} \label{section6}
This paper investigated the generalization performance of regularized least square regression with Fredholm kernel. Generalization bound  is presented for the Fredholm learning model, which shows that the fast learning rate with $O(l^{-1})$ can be reached. In the future, it is interesting to investigate the leaning performance of ranking \cite{hong2} with Fredholm kernel.

\subsection*{Acknowledgments}
The authors would like to thank Prof.Dr.L.Q. Li for his valuable suggestions. This work was supported by the National Natural Science Foundation of China(Grant Nos. 11671161) and the Fundamental Research Funds for the Central Universities (Program Nos. 2662015PY046, 2014PY025).

\end{document}